\documentclass[12pt]{article}
\newtheorem{thm}{Theorem}[section]
\newtheorem{lem}[thm]{Lemma}
\newtheorem{cor}[thm]{Corollary}

\begin{document}
\setcounter{page}{1}

\vspace{2cm}

{\centerline{\bf REMARKS ON GENERALIZED HARDY ALGEBRAS}}

\vspace{10mm}

{\centerline{ROMEO ME\v{S}TROVI\'C${}^{*}$, \v{Z}ARKO PAVI\'CEVI\'C${}^{**}$}}

{\centerline{ AND NOVO LABUDOVI\'C${}^{***}$}}

{\renewcommand{\thefootnote}{}\footnote{1991 {\it 
Mathematics Subject Classification} Primary 46E30. Secondary 30H05, 46J15.}
\setcounter{footnote}{0}}

      \begin{abstract}
For a measure space $(\Omega, \Sigma, \mu)$ with a positive finite measure
$\mu$, and a positive real number $p$, we define the space 
$L_p^{+}(\mu)=L_p^{+}$ of 
all (equivalence classes of) $\Sigma$-measurable complex functions
$f$ defined on $\Omega$ such that the function $\left(\log^+|f|\right)^p$ is 
integrable with respect to $\mu $. 
 We define the metric $d_p$ on $L^{+}_p$
which generalizes the metric introduced by Gamelin and Lumer in [G]
for the case $p=1$. It is shown that the space $L^{+}_p$ is a topological
algebra. On the other hand,  one can define on the space $L_p^{+}$
an equivalent $F$-norm $\vert \cdot\vert_p$ that  makes $L_p^{+}$ into
an Orlicz space. For the case of the normalized Lebesgue's measure 
$dt/2\pi$ on $[0,2\pi)$, it follows that the class 
$N^p(1<p<\infty)$ introduced by I. I. Privalov in [P], may be considered 
as a generalization of the Smirnov class $N^+$. Furthermore, $N^p(1<p<\infty)$
with the associated modular  becomes an Hardy-Orlicz class.
Finally, for a strictly positive and measurable on $[0,2\pi)$ function
$w$, we define the generalized Orlicz space $L_p^{w}(\mathrm{d}t/2\pi)=L^w_p$ 
with the modular $\rho^w_p$ given by the function
$\psi_w(t,u)=\big(\log(1+uw(t))\big)^p$, with a "weight" $w$.  
We observe that the space $L^w_p$ is a generalized Orlicz space with respect 
to the modular $\rho^w_p$.  We examine and compare different 
topologies induced on $L^w_p$ by corresponding "weights" $w$.
   
\end{abstract}

\section{Introduction}

For a  measure space $(\Omega, \Sigma, \mu)$ with
 a nonnegative finite, complete measure $\mu$ not vanishing identically,
denote by $L^p(\mu)=L^p$ $(0<p\le\infty)$ the familiar Lebesgue spaces on
$\Omega$. In Section 2, for $p>0$, we define the class ${L}^+_p(\mu)={L}^+_p$ 
of all (equivalence classes of) $\Sigma$-measurable complex functions
$f$ defined on $\Omega$ such that the function $\log^+|f|$
is in $L^p$.
Every space ${L}^+_p$ is an algebra.
 For each $p>0$, in Section 2 we define the metric $d_p$ on $L^+_p$ by
\begin{eqnarray*}
d_p(f,g)&=&\inf_{t>0}\left[t+\mu\left(\{x\in\Omega
:\,|f(x)-g(x)|\ge t\}\right)\right]\\&&+\int_{\Omega}\left|\left(\log^+|f(x)|
\right)^p-\left(\log^+|g(x)|\right)^p\right|\,\mathrm{d}\mu.
\end{eqnarray*}   
The space $L^+_1$ was introduced in [G, p.~122], with the notation $L(\mu)$ 
in [G]. It was proved in [G] that
$L^+_1$ with the topology given by the metric $d_1$ is a topological algebra.
In Section 2 we prove the same statement for every space $L^+_p, p>0$.
We also define two metrics $\rho_p$ and $\delta_p$ on ${L}^+_p$, and we show
that they induce the same topology on ${L}^+_p$ as the initial metric $d_p$.

By analogy with the Hardy algebra $H(\mu)$ defined in [G], in Section 3 we
define the algebra $H_p(\mu)$ with $p>0$.
It is known (see [G]) that the Smirnov class $N^+$ on the unit disk 
$D: |z|<1$ in the complex plane may be considered as the Hardy algebra 
$H(\mathrm{d}\theta/2\pi)$.
The analoguous results are obtained in Section 5 for the algebra $N^p$, $p>1$,
introduced by I.~I.~Privalov with the notation $A_q$ in [P].

In Section 4 we note that the function $\psi :[0,\infty)\mapsto [0,\infty)$
defined as $\psi(t)=\big(\log(1+t)\big)^p$, 
is an {\it Orlicz function} or a $\varphi$-{\it function}. Further, we
observe that the space
 $L^+_p(\mathrm{d}t/2\pi)=L_p$, $p>0$, consisting of all 
complex-valued functions $f$, defined and measurable on $[0,2\pi)$ for which
$$
\big(\|f\|_p\big)^p:=\int_{0}^{2\pi}\left(\log(1+|f(t)|)\right)^p\,
\frac{\mathrm{d}t}{2\pi}<+\infty.\eqno(1.1)
$$
is the {\it Orlicz class} coinciding with the associated 
 {\it Orlicz space} (see [Mu, Definition 1.4, p.~2]),  
whose generalization we give in Section 6. We prove that the  modular convergence $\|\cdot\|_p$ and 
the norm convergence $|\cdot|_p$ are equivalent.
As an application, we show that
there does not exist a nontrivial continuous linear functional on the space
$(L_p,\|\cdot\|_p)$.

For $p>1$, following I.~I.~Privalov (see [P, p.~93]), a function $f$ 
holomorphic in $D$, belongs to the class $N^p$, if there holds
 $$
\sup_{0\leq r<1}\int_0^{2\pi}\left(\log^{+}\big\vert
f\big(re^{i\theta}\big)\big\vert\right)^p\,\frac{\mathrm{d}\theta}{2\pi}
<\infty,
 $$
In Section 5 we note that the algebra $N^p$ may be considered as the 
{\it Hardy-Orlicz space} with the Orlicz function $\psi$ defined in Section 4.
 Identifying a function $f\in N^p$ with its {\it boundary function} $f^*$, 
the space $N^p$ is identical with the closure of the space of all functions
holomorphic in $D$ and continuous in $\bar{D}:|z|\le 1$ in the space
$(L_p(dt/2\pi)\cap N,|\cdot|_p)$. On the other hand, $(N^p,d_p)$ coincides 
with the space $H_p(d\theta/2\pi)$ defined in Section 3, i.e. with the 
closure of {\it the disk algebra} $P(\bar{D})$ in the space 
$(L_p(dt/2\pi), d_p)$.
Therefore, the space $N^p$ may be considered as generalized Hardy algebra.
From this fact, it is easy to show that $N^p$ is an $F$-algebra 
with respect to the $F$-norm $\|\cdot\|_p$ given by 1.3.

In the last section we note that the real function $\psi_w$ with
$p>0$, defined  on  $\Omega \times [0,\infty)$  by the formula
\[
\psi_w(t,u)=\big(\log(1+uw(t))\big)^p,
\]
is a {\it Musielak-Orlicz function}.
For the Lebesgue measure space 

\noindent $(\Omega, [0,2\pi),dt/2\pi)$,
we denote by $L^w_p(\mathrm{d}t/2\pi)=L^w_p$, $p>0$, the class of all 
(equivalence classes of) complex-valued functions $f$, defined and measurable 
on $[0,2\pi)$ for which
\[
\big(\|f\|^w_p\big)^p:=\int_{0}^{2\pi}\left(\log(1+|f(t)|w(t))\right)^p\,
\frac{\mathrm{d}t}{2\pi}<+\infty.
\]
Then $L^w_p$ is {\it the generalized Orlicz class}  with  the modular 
$\|\cdot\|^w_p$, and $L^w_p$ coincides with the
associated  {\it generalized Orlicz space}.
$\|\cdot\|^w_p$ is a modular in the sense of Definition 1.1 in [Mu, p.~1],  
and by $\rho^w_p(f,g)=\big(\|f-g\|^w_p\big)^{\min\{p,1\}}, f,g\in L^w_p$,
is defined an invariant metric on $L^w_p$. By [Mu, p.~2, Theorem 1.5, and 
p.~35, Theorem 7.7], it follows that the functional $|\cdot|^w_p$ 
defined as
 \[
|f|^w_p=\inf\left\{\varepsilon>0:\,\int_{0}^{2\pi}\left(\log\left(1+
\frac{|f(t)|w(t)}
{\varepsilon}\right)\right)^p\,\frac{\mathrm{d}t}{2\pi}\le\varepsilon
\right\},\quad f\in L^w_p,
  \]          
is a complete $F$-norm.  
We prove that  $(L^w_p,\rho^w_p)$ is an $F$-space for any  weight $w$. 

For two weights $w$ and $\omega$  such that  $\log^+\left(w/{\omega}\right)
\in L^p$, we show that $L^{\omega}_p\subset L^w_p$, and $L^{\omega}_p= L^w_p$
if and only if $\log \left(w/{\omega}\right)\in L^p$.
Further, we prove that the topology  defined on $L^{\omega}_p$ by the metric 
$\rho^{\omega}_p$ is stronger than that induced on $L^{\omega}_p$ by the 
metric $\rho^w_p$. If $p\ge 1$ and  $\log\left(w/{\omega}\right)\notin L^p$, 
then $L^{\omega}_p$ is a proper subset of $L^w_p$, and the topology defined on 
$L^{\omega}_p$  by the metric $\rho^{\omega}_p$ is strictly stronger than 
that induced on $L^{\omega}_p$ by the metric $\rho^w_p$.
As an application, we show that if $\log w\in L^p$, then  $(L_p,\rho^w_p)$ 
is an $F$-algebra, and the metric topologies  $\rho^w_p$ and $\rho_p$ 
are the same. Finally, if $p\ge 1$ and $w$ is a weight such that
$\log^+w\in L^p$, we give four equivalent necessary and sufficient conditions
for the space $(L_p,\rho^w_p)$ to be an $F$-algebra.

\section{Equivalent metrics on the space $L^+_p$ $(0<p<\infty)$}
Let $(\Omega, \Sigma, \mu)$ be a measure space, i.e. $\Omega$ is a
nonempty set, $\Sigma$ is a $\sigma$-algebra of subsets of $\Omega$ and
$\mu$ is a nonnegative finite, complete measure not vanishing identically.
Denote by $L^p(\mu)=L^p$ $(0<p\le\infty)$ the familiar Lebesgue spaces on
$\Omega$. For $p>0$, we define the class ${L}^+_p(\mu)={L}^+_p$ of 
all (equivalence classes of) $\Sigma$-measurable complex functions
$f$ defined on $\Omega$ such that the function $\log^+|f|$
is in $L^p$, where $\log^+a=\max\{\log a,0\}$, i.e. such that
 $$
\int_{\Omega}\left(\log^+|f(x)|\right)^p\,
\mathrm{d}\mu<+\infty.
 $$
 Obviously, $L^+_q\subset L^+_p$ for $q>p$,
and from the inequality $\log^+x\le x^s/se, x\le 0, s>0$, we see that 
$\bigcup_{p>0}L^p\subset\bigcap_{p>0}L^+_p$. Combining the inequalities
$\log^+|f+g|\le\log^+|f|+\log^+|g|+\log 2$, $\log^+|fg|\le\log^+|f|+\log^+
|g|$ with $(|x|+|y|)^p\le 2^{\max\{p-1,0\}}(|x|^p+|y|^p)$ and
$(|x|+|y|+|z|)^p\le 3^{\max\{p-1,0\}}(|x|^p+|y|^p+|z|^p)$, respectively,
we see that every space ${L}^+_p$ is an algebra with respect to the pointwise
addition and multiplication.
 For each $p>0$, we define the metric $d_p$ on $L^+_p$ by
\begin{eqnarray*}
d_p(f,g)&=&\inf_{t>0}\left[t+\mu\left(\{x\in\Omega
:\,|f(x)-g(x)|\ge t\}\right)\right]\\
\qquad \qquad&&+\int_{\Omega}\left|\left(\log^+|f(x)|
\right)^p-\left(\log^+|g(x)|\right)^p\right|\,\mathrm{d}\mu.\qquad
\qquad {\rm(2.1)}
\end{eqnarray*} 
The space $L^+_1$ was introduced in [G, p.~122], with the notation $L(\mu)$ 
in [G]. In fact, the above metric $d_p$ with $p=1$ coincides with the
Gamelin-Lumer's  metric $d$ defined on $L^+_1$. It was proved in 
 [G, p.~122, Theorem 2.3] that the space
$L^+_1$ with the topology given by the metric $d_1$ is a topological algebra.
The following result is a generalization of the corresponding result
for the case $p=1$. The proof of this result is completely analogous to
those for the case $p=1$ given in [G, p.~122]{3}, 
and therefore, may be omitted.

\begin{thm}
The space $L^+_p$ with the metric $d_p$ given by (1.2) is a {\it topological 
algebra}, i.e. a topological vector space with a complete metric in which
multiplication is continuous.  
\end{thm}
By the inequality 
 $$
\big(\log(1+|x|)\big)^p\le 2^{\max\{p-1,0\}}\big((\log 2)^p+(\log^+|x|)^p\big),
\eqno(2.2)
 $$ it follows that $f$ belongs to ${L}^+_p$
if and only if there holds
 $$
\big(\|f\|_p\big)^p:=\int_{\Omega}\left(\log(1+|f(x)|)\right)^p\,\mathrm{d}\mu
<\infty.
 $$
By the inequality $\log(1+|f+g|)\le\log(1+|f|)+(\log(1+|g|)$ and Minkowski's
inequality, it follows that for $p>1$ the function $\rho_p$ defined
as
  $$
\rho_p(f,g)=\|f-g\|_p,\quad f,g\in L^+_p,\eqno(2.3)
 $$
satisfies the triangle inequality. Combining the above inequality with 
$\big(|x|+|y|\big)^p\le |x|^p+|y|^p$ for $0<p\le 1$, we see that  
the function $\rho_p$ defined  as
 $$
\rho_p(f,g)=\big(\|f-g\|_p\big)^p,\quad f,g\in L^+_p, 0<p\le 1,
 $$
satisfies also the triangle inequality. Hence,  $\rho_p$ is
an invariant metric on $L^+_p$ for all $p>0$.

Recall that a subset $K$ of $L^p$ forms an
{\it uniformly integrable family} if for given $\varepsilon> 0$ there exists a
$\delta> 0$ so that
\[
\int_E\vert f(x)\vert\,\mathrm{d}\mu<\delta
\quad\mbox{for all }\quad f\in K,
\]
whenever $E\subset\Omega$ with  its measure $\mu (E)<\delta$.        

 Two metrics (or norms) defined on the same space will be called {\it  
equivalent} if they induce the same topology.

For the proof of Theorem 2.3 we will need the following lemma. 
\begin{lem}
{\rm ([G, p.~122, Theorem 1.3]).} Let $(f_n)$ be a sequence in $L^1$ and 
$f\in L^1$ such that $f_n\to f$ in $L^1$. Then a sequence $(f_n)$ is a 
uniformly integrable family. Conversely, if a sequence $(f_n)$ is a uniformly
integrable family on $\Omega$, and $f_n\to f$ in measure, then $f$ belongs
to $L^1$ and $f_n\to f$ in $L^1$.
\end{lem}
\begin{thm}
The metric $d_p$  defines a topology for $L^+_p$ which is
equivalent to the topology defined by the metric $\rho_p$.
\end{thm}
{\it Proof.}
 Suppose first that $\rho_p(f_n,f)\to 0$ as $n\to\infty$, where $(f_n)$
and $f$ are in $L^+_p$. Then by Chebyshev's inequality, it is easily seen 
that $f_n\to f$ in measure on $\Omega$. For simplicity, put 
\[
\big(\| f\|_E\big)^{\max\{p,1\}}=\int_E\left(\log(1+|f(x)|)\right)^p
\,\mathrm{d}\mu
\]
for any measurable set $E\subset\Omega$. 
By the triangle inequality and (2.2), we have
\begin{eqnarray*}   
\| f_n\|_E &\le & \|f_n-f\|_E+\| f\|_E\\
    \qquad\qquad\qquad       &\le& \rho_p(f_n,f)+2^{\max\{1-1/p,0\}}\qquad\qquad
\qquad\qquad\qquad {\rm (2.4)}\\
\qquad&\times&\left(\mu(E)(\log 2)^p+
\int_E\left(\log^+|f(x)|\right)^p\,\mathrm{d}\mu\big)^{1/\max\{p,1\}}\right).
    \end{eqnarray*}   
This shows that $\big\{\left(\log^+|f_n(x)|\right)^p:\, n\in {\bf{\rm N}}
\big\}$
 form a uniformly integrable family, and by Lemma 2.2, $d_p(f_n,f)\to 0$ 
as $n\to\infty$. 

Conversely, assume that $d_p(f_n,f)\to 0$ as $n\to\infty$. Then, by the
definition of  $d_p$, $f_n\to f$ in measure, and by Lemma 2.2,  
$\big\{\left(\log^+|f_n(x)|\right)^p:\, n\in {\bf{\rm N}}\big\}$
are uniformly integrable. Replacing in (2.3) $f$ by $f_n$ and $f_n$
by $f-f_n$, we see that the family 
$\big\{\left(\log(1+|f_n(x)-f(x)|)\right)^p:\, n\in {\bf{\rm N}}\big\}$
is uniformly integrable. Thus, by Lemma 2.2, 
$\rho_p(f_n,f)\to 0$ as $n\to\infty$. Hence the metrics $d_p$ and $\rho_p$
are equivalent.\\

\noindent{\it Remark.} {\it Using the same argument applied in the proof of
Theorem 2.3, it is easy to see that the metrics $\rho_p$ and $d_p$
are equivalent with the metric $\delta_p$ given on $L^+_p$ by}
\begin{eqnarray*}
\delta_p(f,g)&=&\inf_{t>0}\left[t+\mu\left(\{x\in\Omega
:\,|f(x)-g(x)|\ge t\}\right)\right]\\
&&+\left(\int_{\Omega}\left|\log^+|f(x)|
-\log^+|g(x)|\right|^p\,\mathrm{d}\mu\right)^{1/\max\{p,1\}},
f,g\in L^+_p.
\end{eqnarray*}  \\

\noindent{\it Remark.} {\it In} [Y, Remark  5, p.~{\rm 460], 
{\it M.~Hasumi pointed out
that the Yanagihara's metric $\rho=\rho_1$  defines a topology for 
the space $L^+_1=L(\mu)$, which is equivalent to the metric  topology  $d_1=d$
 used by Gamelin-Lumer in} [G, p.~122].   

\begin{cor} The space $L^+_p$, $p>0$, with the topology given by the
metric $\rho_p$ is an $F$-algebra, i.e. a topological algebra with a 
complete translation invariant metric $\rho_p$.  
\end{cor}
{\it Proof.}
We see from Theorems 2.1 and 2.3 that it is sufficient to show that the space 
$L^+_p$ is complete with respect to the metric $\rho_p$. The completeness of
$L^+_p$  may be proved by the standard manner as the 
completeness of the Lebesgue spaces $L^p$ or an arbitrary generalized Orlicz 
space with a corresponding $F$-norm (for example, see the proof of Theorem
7.7 in [Mu, p.~35]).
In view of this, note that
$L^+_p$ may be considered as the generalized Orlicz space  $L^w_p$
with a constant function $w(t)\equiv 1$ on $[0,2\pi)$ (see Section 6).

\begin{thm}
If $0<p<s<\infty$, then $L^+_s\subset L^+_p$, and the metric 
$d_p$ induces a topology for $L^+_s$ which is
coarser than the topology defined on $L^+_s$ by the initial metric $d_s$.
\end{thm}
{\it Proof.}
If $0<p<s<\infty$, $L^+_s\subset L^+_p$ is obvious.
Let $(f_n)$ be a sequence in $L^+_s$  and 
$f\in L^+_s$ such that $d_s(f_n, f)\to 0$ as $n\to 0$.
This means by definition of $d_s$, that  $\left(\log^+|f_n|\right)^s
\to \left(\log^+|f|\right)^s$ in $L^1$, and  $f_n\to f$ in measure.
Then by Lemma 2.2, a sequence $\left(\left(\log^+|f_n|\right)^s\right)$ is a 
uniformly integrable family on $\Omega$. Since $0<p<s<\infty$,
it is routine to verify that $\left(\left(\log^+|f_n|\right)^p\right)$ is also 
a uniformly integrable family on $\Omega$.
Then by Lemma 2.2, $\left(\log^+|f_n|\right)^p\to \left(\log^+|f|\right)^p$ in
 $L^1$, and hence $d_p(f_n, f)\to 0$ as $n\to 0$.


\section{A generalization of Hardy algebras}

Let $A$ be a {\it uniform algebra} on a compact Hausdorff space $\Omega$.
Let $\mu$ be a {\it representing measure} for a nonzero complex-valued 
homomorphism  $\alpha$ of $A$. For $0<q\le\infty$, denote by $H^q(\mu)$ 
{\it the Hardy space}
(the closure of $A$ in the Lebesgue space $L^q=L^q(\mu))$. Motivated
by the definition of the {\it Hardy algebra} $H(\mu)$ given in 
[G, p.~123], 
here we define the space $H_p(\mu)=H_p(0<p<\infty)$, as the $d_p$-closure
of $A$ in the space $L^+_p(\mu)=L^+_p$ defined in the Section 2.
Observe that the space $H_1$ coincides with the Hardy algebra
$H(\mu)$ defined in [G, p.~122]. From Theorem 2.5 we obtain immediately its 
$H_p$-analogue.
\begin{thm}
If $0<p<s<\infty$, then $H_s\subset H_p$,  and the metric $d_p$ induces a 
topology for $H_s$ which is coarser than the topology defined on $H_s$ by the 
initial metric $d_s$.
\end{thm}
Denote by $M_{\alpha}$ the set of all representing measures for a
homomorphism $\alpha$. For the proof of Theorem 3.2, we will need the
following lemma which generalizes Corollary 2.2 in [G, p.~122]. 
\begin{lem} For all $0< q\le\infty$ and all $p>0$, the topology of $L^q$
is stronger than the topology of the space $L^+_p$.
\end{lem}
{\it Proof.}
Clearly, it suffices to consider the case $0< q\le 1$. 
Let $(f_n)$ be a sequence in $L^q$ such that $f_n\to f$ in $L^q$
for some $f\in L^q$. From the inequality $\log(1+x)\le x^s/s$,
$x\ge 0, 0<s\le 1$, it follows that  $\big(\log(1+x)\big)^p\le px^q/q$,
and hence $\rho_p(f_n,f)\to 0$ as $n\to \infty$. Therefore, 
we conclude by Theorem 2.3 that $d_p(f_n,f)\to 0$ as $n\to\infty$, 
which completes the proof.

\begin{thm}
Suppose that the space $M_{\alpha}$ is finite-dimensional, and let
$\mu$ be a  core representing measure for a  homomorphism
 $\alpha$. Then for any fixed $0< q\le\infty$, there holds
 $$
\bigcap_{p>0} H_p\bigcap L^q=H^q(\mu).\eqno(3.1)
 $$
Furthermore, there holds
 $$
 H_p\bigcap L^q=H^q(\mu)\quad\mbox{for }1\le p\le q\le\infty.\eqno(3.2)
 $$
\end{thm}
{\it Proof.}
By [G, p.~125, Theorem 4.2], for all $0<q\le\infty$ we have 
$H_1\bigcap L^q=H^q(\mu)$. This shows that $\bigcap_{p>0} H_p\bigcap L^q
\subseteq H^q(\mu)$. 

Suppose now that $f\in H^q(\mu)$ for some $0<q\le\infty$. Then there
is a sequence $(f_n)$ in  $A$ such that $f_n\to f$ in $L^q$ as $n\to\infty$.
By Lemma 3.2,  $f_n\to f$ in $L^+_p$ as $n\to\infty$, for each $p>0$.
Thus, $f\in \bigcap_{p>0} H_p$,  which implies that
$H^q(\mu)\subseteq\bigcap_{p>0} H_p\bigcap L^q$.
This yields (3.1). If $1\le p\le q\le\infty$, by Theorem (3.1) we have 
 $$
H_p\bigcap L^q\subseteq H_1\bigcap L^q =H^q(\mu).
$$
Conversely, from (3.1) we obtain
$H^q(\mu)=\bigcap_{s>0} H_s\bigcap L^q\subseteq H_p\bigcap L^q$.
This proves (3.2).

\section{Orlicz spaces $L^+_p(\mathrm{d}t/2\pi)(0<p<\infty)$}
The function $\psi :[0,\infty)\mapsto [0,\infty)$
defined as $\psi(t)=\big(\log(1+t)\big)^p$,
is continuous and nondecreasing in $[0,\infty)$, such that $\psi(0)=0$, 
$\psi(t)>0$ for $t>0$, and 
$\lim_{t\to +\infty}\psi(t)=+\infty$, is called an {\it Orlicz function}
or a $\varphi$-{\it function} (see  [Mu, p.~4, Examples 1.9]). Further,
observe that the space
 $L^+_p(\mathrm{d}t/2\pi)=L_p$, $p>0$, consisting of all 
complex-valued functions $f$, defined and measurable on $[0,2\pi)$, for which
 $$
\big(\|f\|_p\big)^p:=\int_{0}^{2\pi}\left(\log(1+|f(t)|)\right)^p\,
\frac{\mathrm{d}t}{2\pi}<+\infty.
 $$
is the {\it Orlicz class} (see [Mu, p.~5]), whose generalization we
give in Section 6.
It follows by the dominated convergence theorem
that the class $L_p$ coincides with the
associated  {\it Orlicz space} (see [Mu, p.~2, Definition 1.4]), consisting 
of those functions $f\in L_p$ such that
\[
\int_{0}^{2\pi}\left(\log(1+c|f(t)|)\right)^p\,
\frac{\mathrm{d}t}{2\pi}\to 0\quad {\mbox as }\quad c\to 0+.
\]
 Furthermore, 
since by $\rho_p(f,g)=\big(\|f-g\|_p\big)^{\min\{p,1\}}, f,g\in L_p$,
is defined an invariant metric on $L_p$, the function
$\|\cdot\|_p$ given by 4.1
is a {\it modular} in the sense of Definition 1.1 in [Mu, p.~1]. 
For any fixed $f\in L_p$, by the monotone convergence theorem, we see that 
$\lim_{c\to 0}\psi(cf)\to 0$, and so $(L_p,\rho_p)$ is a {\it modular
space} in the sense of Definition 1.4 in [Mu, p.~2].
 In other words,
the function $\|\cdot\|_p$  is an $F$-norm. It is well known (see [Mu, 
Theorem 1.5, p.~2 and Theorem 7.7, p.~35]) that the functional $|\cdot|_p$ 
defined as
 \[
|f|_p=\inf\left\{\varepsilon>0:\,\int_{0}^{2\pi}\left(\log\left(1+\frac{|f(t)|}
{\varepsilon}\right)\right)^p\,\frac{\mathrm{d}t}{2\pi}\le\varepsilon
\right\},\quad f\in L_p,
  \]          
is a complete $F$-{\it norm}. Furthermore (see [L, p.~54]), $L_p$ is a 
completion (closure) of the space of all continuous functions on $[0,2\pi)$
in the space $(L_p,|\cdot|_p)$.
\begin{thm} The $F$-norms $\|\cdot\|_p$ and $|\cdot|_p$ induce the same 
topology on the space $L_p$. In other words, the norm and modular convergences
are equivalent.
\end{thm} 
{\it Proof.}
Since a $\varphi$-function $\psi(t)=\big(\log(1+t)\big)^p$  satisfies
obviously the so-called $\Delta_2$-{\it condition} given by
$$
\phi(2t)\le c\phi(t)\eqno (\Delta_2)
$$
with a constant $c=2^p$, the assertion follows immediately from
[L, p.~55, 2.4]. But, we give here a direct proof.

Let $K_{\rho}(0,r)=\{f\in L_p:\, \|f\|_p<r\}$ be a neighborhood of $0$ in the 
space $(L_p,\|\cdot\|_p)$, and let $f\in K_{\rho}(0,r)$ be any fixed.
  Consider a neighborhood
 $K_{\psi}(f,\varepsilon/2)=\{f\in L_p:\, \|f\|_p<r\}$ of $f$ in the space
$(L_p,|\cdot|_p)$, and assume $g\in K_{\psi}(f,\varepsilon/2)$. Then we
have
\[
\frac{\varepsilon}{2}>|g-f|_p=\inf\big\{\delta >0:\,\left\|\frac{g-f}{\delta}
 \right\|<\delta\big\}=\sigma. 
\]
Thus, from  $\varepsilon/2>\sigma$ we infer that
\[
\left\|\frac{g-f}{\varepsilon/2}
 \right\|<\frac{\varepsilon}{2},
\]
whence we obtain
\begin{eqnarray*}
\|g\|_p&\le& \|g-f\|_p+\|f\|_p\\
&\le& \left\|\frac{g-f}{\varepsilon/2}
 \right\|+\left|\frac{\varepsilon}{2}\right|+r-\varepsilon\\
&<& \frac{\varepsilon}{2}+ \log\big(1+\frac{\varepsilon}{2}\big)+r
-\varepsilon\\
&<& \frac{\varepsilon}{2}+ \frac{\varepsilon}{2}+r-\varepsilon=r.
\end{eqnarray*}
This shows that $K_{\psi}(f,\varepsilon/2)\subset K_{\rho}(0,r)$. 

Conversely, let $K_{\psi}(0,r)$  be a neighborhood of $0$ in the 
space $(L_p,|\cdot|_p)$, and let $f\in K_{\psi}(0,r)$ be any fixed.
Put $\sigma =r-|f|_p$, take $k$ an integer  so that $k> 2/\sigma$,
and set $\varepsilon=\sigma/2k$. Then
for all $g\in K_{\rho}(f,\varepsilon)$ by the triangle inequality, we have
\[
\left\|\frac{g-f}{\sigma/2}\right\|<\|k(g-f)\|_p\le k\|g-f\|_p\le 
k\varepsilon=\frac{\sigma}{2}.
\]
Hence, we obtain  $|g-f|_p\le \sigma/2$, and so
\[
|g|_p\le |g-f|_p+|f|_p\le\frac{\sigma}{2}+r-\sigma=r-\frac{\sigma}{2}<r. 
\]
Thus, $g\in K_{\psi}(0,r)$, i.e. $K_{\rho}(f,\varepsilon)\subset 
 K_{\psi}(0,r)$. Hence, the topologies defined on the space $L_p$ 
by $\|\cdot\|_p$ and $|\cdot|_p$ are the same.\\ 
   
As an application of Theorem 4.1, we obtain the following result.
\begin{cor}
There does not exist a nontrivial continuous linear functional on the space
$(L_p,\|\cdot\|_p)$.
\end{cor}
{\it Proof.}
Since the Orlicz function $\psi(t)=\big(\log(1+t)\big)^p$ satisfies
the condition $\liminf_{t\to +\infty} t^{-1}\psi(t)=0$, it
follows by [MO] that there does not exist a nontrivial modular 
continuous linear functional on the space $(L_p,|\cdot|_p)$.  
Hence, this is true by Theorem 4.1 for the space $(L_p,\|\cdot\|_p)$.

\section{Privalov's spaces $N^p(1<p<\infty)$}
As in [G, p.~125], {\it the Hardy algebra} $H_1(d\theta/2\pi)$
consists of all functions $f$ holomorphic on the unit disk $D:|z|<1$
for which the family
\[
\left\{\log^+\vert f(re^{i\theta})\vert:\,0\leq r<1\right\}
\] is uniformly integrable on the unit circle $T$, that is, for a given 
$\varepsilon >0$, there exists $\delta>0$ so that
\[
\int_E\log^+\big\vert f\big(re^{i\theta}\big)\big\vert\,\frac{\mathrm{d}
\theta}{2\pi}<\varepsilon, \quad 0\leq r<1,
\]
for any measurable set $E\subset [0,2\pi)$ with its Lebesgue measure 
$|E|<\delta$. This space is usually called {\it the Smirnov class} $N^+$.
For $p>1$, following I.~I.~Privalov (see [P, p.~93], where $N^p$ is
denoted as $A_q$), a 
function $f$ holomorphic in $D$, belongs to the class $N^p$, if there holds
 $$
\sup_{0\leq r<1}\int_0^{2\pi}\left(\log^{+}\big\vert
f\big(re^{i\theta}\big)\big\vert\right)^p\,\frac{\mathrm{d}\theta}{2\pi}
<\infty,
 $$
For $p=1$, the condition (5.1) defines  {\it the Nevanlinna class} $N$ of
holomorphic functions in $D$. Recall that for $f\in N$, {\it the radial limit}
 $$
 f^{*}(e^{i\theta})=\lim_{r\to 1} f\big(re^{i\theta}\big)
 $$
exists for almost every $e^{i\theta}$ and $\log\vert f^*\vert\in L^1$
unless $f\not\equiv 0$. It is known (see [Mo]) that
 \[
\bigcup_{p>1}N^p\subset N^+\subset N.
\]
Observe that the algebra $N^p$ may be considered as the {\it Hardy-Orlicz
space} with the Orlicz function $\psi :[0,\infty)\mapsto [0,\infty)$
defined as $\psi(t)=\big(\log(1+t)\big)^p$. For more informations on the
Hardy-Orlicz spaces see [Mu, Ch. IV, Sec. 20]. Identifying a function
$f\in N$ with its {\it boundary function} $f^*$, by [L, 3.4, p.~57], the 
space $N^p$ is identical with the closure of the space of all functions
holomorphic in $D$ and continuous in $\bar{D}:|z|\le 1$ in the space
$(L_p(dt/2\pi)\cap N,|\cdot|_p)$. On the other hand, $(N^p,d_p)$ coincides 
with the space $H_p(d\theta/2\pi)$ defined in Section 3, i.e. with the 
closure of 
{\it the disk algebra} $P(\bar{D})$ in the space $(L_p(dt/2\pi), d_p)$.
Therefore, the space $N^p$ may be considered as generalized Hardy algebra.
From this fact, Theorem 2.3 and Corollary 2.4, it follows immeditately
the following result of M.~Stoll obtained in [S].
\begin{thm}
{\rm ([S, Theorem 4.2]).} For  $p>1$, the space $N^p$ with the topology
given by the metric $\rho_p$ is an $F$-algebra. Furthermore, the 
polynomials are dense in $(N^p,\rho_p)$, and hence $N^p$ is separable.
\end{thm}
\section{The generalized Orlicz spaces $L_p^w$}
Let $(\Omega, \Sigma, \mu)$ be a measure space, and let $w$ be a 
$\Sigma$-measurable nonnegative function defined on $\Omega$ such that
$w(t)>0$ $\mu$-almost every on $\Omega$. For $p>0$,
the real function $\psi_w$ defined on  $\Omega \times [0,\infty)$
by the formula 
\[
\psi_w(t,u)=\big(\log(1+uw(t))\big)^p
\]
is a {\it Musielak-Orlicz function}, since it satisfies {\it the Caratheodory
conditions}, i.e. it satisfies the following conditions (see [Mu, p.~33, 
Definition 7.1]):

(i) $u\mapsto \psi_w(t,u)$ is a $\varphi$-function of the variable $u\ge 0$
for $\mu$-almost every $t$, i.e. is a nondecreasing, continuous function of
$u$ such that $\psi_w(t,0)=0$, $\psi_w(t,u)>0$ for $u>0$,
     $\lim_{u\to\infty}\psi_w(t,u)\to\infty$;

(ii) $\psi_w(t,u)$ is a $\Sigma$-measurable function of $t$ for all
$u\ge 0$.

For the Lebesgue's measure space $(\Omega, [0,2\pi),dt/2\pi)$,
denote by $L^w_p(\mathrm{d}t/2\pi)=L^w_p$, $p>0$, the class of all 
(equivalence classes of) complex-valued functions $f$, defined and measurable 
on $[0,2\pi)$ for which
\[
\big(\|f\|^w_p\big)^p:=\int_{0}^{2\pi}\left(\log(1+|f(t)|w(t))\right)^p\,
\frac{\mathrm{d}t}{2\pi}<+\infty.
\]
$L^w_p$ is called {\it the generalized Orlicz class}  with  the modular 
$\|\cdot\|^w_p$ (see [Mu,  p.~33, Definition 7.2]). It follows by the 
dominated convergence theorem that the class $L^w_p$ coincides with the
associated  {\it generalized Orlicz space} consisting of those functions
$f\in L^w_p$ for which
\[
\int_{0}^{2\pi}\left(\log(1+c|f(t)|w(t))\right)^p\,
\frac{\mathrm{d}t}{2\pi}\to 0\quad\mbox{as }\quad c\to 0+.
\]
Furthermore, from the inequality
$\|cf\|^w_p\le\|c\|^w_p+\|f\|^w_p$, we see that $L^w_p$ coincides with
 the space of all {\it finite elements} of $L^w_p$ consisting
of those functions $f\in L^w_p$ such that $cf\in L^w_p$ for every $c>0$.
Since $w(t)>0$ almost every on $[0,2\pi)$ and
\[
\|f+g\|^w_p=|(f+g)w|_p\le |fw|_p+|gw|_p=\|f\|^w_p+\|g\|^w_p,
\]
$\|\cdot\|^w_p$ 
is a  modular in the sense of Definition 1.1 of [Mu, p.~1]. 
In other words, the function $\|\cdot\|^w_p$  is an $F$-norm.
 Hence, 
 by $\rho^w_p(f,g)=\big(\|f-g\|^w_p\big)^{\min\{p,1\}}, f,g\in L^w_p$,
is defined an invariant metric on $L^w_p$. By [Mu, p.~2, Theorem 1.5, and 
p.~35, Theorem 7.7], it follows that the functional $|\cdot|^w_p$ 
defined as
 \[
|f|^w_p=\inf\left\{\varepsilon>0:\,\int_{0}^{2\pi}\left(\log\left(1+
\frac{|f(t)|w(t)}
{\varepsilon}\right)\right)^p\,\frac{\mathrm{d}t}{2\pi}\le\varepsilon
\right\},\quad f\in L^w_p,
  \]          
is a complete $F$-norm.  
\begin{thm}
$(L^w_p,|\cdot|^w_p)$ is an $F$-space. Furthermore, $F$-norms $\|\cdot\|^w_p$ 
and $|\cdot|^w_p$ induce the same topology on the space $L^w_p$.
\end{thm}  
{\it Proof.}
 Since the functional $|\cdot|_p$ is a complete $F$-norm, and hence
for a fixed $f\in L^w_p$, $c\mapsto cf$
is a continuous mapping from ${\bf{\rm C}}$ into $L^w_p$, it suffices
(see [DS, p.~51])  to check the continuity of the mapping
 $f\mapsto cf$ from $L^w_p$ into $L^w_p$, for a fixed $c\in{\bf{\rm C}}$.
Take $k$ an integer with $|c|\le k$. Then
    \[
|cf|\le |kf|\le k|f|;
   \]   
 so $f\mapsto cf$ is continuous.

The second assertion of the theorem can be proved completely analogously
as Theorem 4.1.

\begin{thm}  For any  weight $w$, $(L^w_p,\rho^w_p)$ is an $F$-space. 
\end{thm}
{\it Proof.}
To show that $(L^w_p,\rho^w_p)$ is an $F$-space, 
it suffices by Theorem 6.1 to prove that $(L^w_p,\rho^w_p)$ is complete.
Let $(f_n)$ be a Cauchy sequence in $(L^w_p,\rho^w_p)$. This means
that $(wf_n)$ is a Cauchy sequence in $(L_p,\rho_p)$, and by Corollary 2.4,
there is a $g\in L_p$ such that $\rho_p(wf_n,g)\to 0$ as $n\to\infty$.
Then $\rho^w_p(f_n,g/w)=\rho_p(wf_n,g)\to 0$, and since $\rho^w_p(g/w,0)=
\rho_p(g,0)<\infty$, we see that $f=g/w\in L^w_p$. Hence, the space
$(L^w_p,\rho^w_p)$ is complete. 
\begin{thm}
Let  $w$ and $\omega$ be two weights such that  $\log^+\left(w/{\omega}\right)
\in L^p$. Then:

{\rm(i)}  $L^{\omega}_p\subset L^w_p$, and $L^{\omega}_p= L^w_p$
if and only if $\log \left(w/{\omega}\right)\in L^p$.

{\rm(ii)}
The topology  defined on $L^{\omega}_p$ by the metric $\rho^{\omega}_p$
 is stronger than that induced on $L^{\omega}_p$ by the metric $\rho^w_p$.

{\rm(iii)} Let $p\ge 1$. If  $\log\left(w/{\omega}\right)\notin L^p$, then
 $L^{\omega}_p$ is a proper subset of $L^w_p$.
Furthermore, the topology defined on $L^{\omega}_p$  by the metric 
$\rho^{\omega}_p$ is strictly stronger than that induced on $L^{\omega}_p$ by 
the metric $\rho^w_p$.

{\rm(iv)} Let $p\ge 1$. If $(L^{\omega}_p,\rho^w_p)$ with $p\ge 1$ is a 
complete metric space, then $\log \left(w/{\omega}\right)\in L^p$, and 
the topology  defined on $L^{\omega}_p=L^w_p$ by the metric $\rho^{\omega}_p$
coincides with that induced on $L^{\omega}_p$ by the metric $\rho^w_p$.
\end{thm}
{\it Proof.} 
(i)
The inclusion relation  $L^{\omega}_p\subset L^w_p$ follows immediately from 
the inequality $\|f\|^w_p=\|f\omega w/\omega\|_p\le \|f\omega\|_p+\|w/\omega\|
_p\|f\|^{\omega}_p+\|w/\omega\|_p=$, $f\in L^{\omega}_p$, and the fact 
that $w/\omega\in L_p$. If $L^{\omega}_p= L^w_p$, then since $1/w\in 
L^w_p=L^{\omega}_p$, it follows that $\log^+\left({\omega}/w\right)\in L^p$, 
and so $\log\left(w/{\omega}\right)\in L^p$. Conversely, if $\log\left(
w/{\omega}\right)
\in L^p$, then $\log^+\left({\omega}/w\right)\in L^p$, and thus
$L^w_p\subset L^{\omega}_p$, which implies $L^{\omega}_p=L^w_p$.

(ii) Assume that $(f_n)$ is a sequence in $L^{\omega}_p$ and 
$f\in L^{\omega}_p$  such that $f_n\rightarrow f$ in $(L^{\omega}_p,\rho_p)$ as
$n\rightarrow\infty$. This means that $\rho^{\omega}_p(f_n,f)\to 0$,
or equivalently, $\rho_p(\omega f_n,\omega f)\to 0$ as $n\rightarrow\infty$.
Since $w/\omega\in L_p$, by Corollary 2.4, we have $\rho_p(wf_n,wf)\to 0$, 
as $n\rightarrow\infty$. Thus, $\rho^w_p(f_n,f)\to 0$  as
$n\rightarrow\infty$. 

(iii) If $\log\left(w/{\omega}\right)\notin L^p$, the strict inclusion 
relation  $L^{\omega}_p\subset L^w_p$ follows from the assertion (i). Since
$\log^+ \left(w/{\omega}\right)\in L^p$ and $\log\left(w/{\omega}\right)
\notin L^1$, by [MP, Theorem 3.1], there holds
\[
 \inf_{P\in{\cal P}_0}\int_0^{2\pi}\left(\log
\left(1+\vert P(e^{it})\vert \frac{w(t)}{\omega(t)}\right)\right)^p\,
\frac{\mathrm{d}t}{2\pi}=0,
\]
or equivalently, that there exists a sequence $\left(P_n\right)$ of (analytic)
polynomials with $P_n(0)=1$ such that $\rho^w_p\left(P_n\omega,0\right)\to 0$ 
as $n\to\infty$. On the other hand, by the same theorem, there holds
\[    
\inf_{P\in{\cal P}_0}\int_0^{2\pi}\left(\log\left(1+\vert P(e^{it})\vert
\right)\right)^p\,\frac{\mathrm{d}t}{2\pi}=(\log 2)^p,
  \]
 and hence the above sequence $\left(P_n\right)$ does not converge to $0$ in 
the space $(L_p,\rho_p)$, or equivalently that a
sequence $\left(P_n/{\omega}\right)$ does not converge to $0$ in 
the space $(L^{\omega}_p,\rho^{\omega}_p)$ Hence, the metric topology 
$\rho^{\omega}_p$ is strictly stronger than the metric topology $\rho^w_p$. 

(iv) Suppose that $(L^{\omega}_p,\rho^w_p)$ is a complete metric space,
but $\log \left(w/{\omega}\right)\notin L^p$. We will consider two cases.

{\it Case} 1. $\log \left(w/{\omega}\right)\notin L^1$.
 By Theorem 6.2, the space  $(L^{\omega}_p,\rho^{\omega}_p)$ is an $F$-space,
 and by the assumption, $(L^{\omega}_p,\rho^w_p)$ is complete,
and henceforth it is an $F$-space. We see from  (iii) that
the topology defined on $L^{\omega}_p$  by the metric 
$\rho^{\omega}_p$ is strictly stronger than that induced on $L^{\omega}_p$ by 
the metric $\rho^w_p$. Consider the identity map
\[
j:(L^{\omega}_p,\rho^{\omega}_p)\to (L^{\omega}_p,\rho^w_p).
\] 
  Then by (ii), $j$ is continuous. By the open mapping theorem, we conclude 
that the inverse $j^{-1}$
of $j$ is also continuous. This shows that the metric topologies
$\rho^{\omega}_p$ and $\rho^w_p$ must be the same. A contradiction.

{\it Case} 2. $\log \left(w/{\omega}\right)\in L^1$.
Since  $\log^+\left({\omega}/w\right)\in L^p\subset L^1$, it follows that
$\log^+\left({\omega}/w\right)\notin L^1$, and hence $1/w\in 
 L^w_p\setminus L^{\omega}_p$.
 Define the outer function $F$  by 
    \[
F(z)=\exp\left(\int_0^{2\pi}\frac{e^{it}+z}{e^{it}-z}\log\frac{w(t)}
{\omega(t)}\right)\,\frac{dt}{2\pi},\quad z\in D.
    \] 
Then $|F^*(e^{it})|=w(t)/{\omega(t)}$ at almost every $t\in [0,2\pi)$.
If $p=1$, by the canonical factorization theorem for the Smirnov class
 $N^+$ (see [D, Theorem 2.10]), it folows that $F$ belongs to $N^+$.
Similarly, if $p>1$, by the canonical factorization theorem for the Privalov's
spaces $N^p$ (see [P, p.~98]), it folows that $F$ belongs to $N^p$.
For simplicity, we shall write  $N^1$ instead of $N^+$.
By a result of Mochizuki  (see [Mo, Theorem 2]), there is a sequence $(f_n)$ 
in $N^p$, $p>1$, such that $f_nF\rightarrow 1$ in $N^p$ as $n\to\infty$. 
For $p=1$, put $f_n=1/F\in N^1$ for all $n$.
 Hence, in both cases  $\rho_p(f^*_nF^*, 1)\to 0$, or equivalently
$\rho_p(f^*_nw/\omega, 1)\to 0$ as $n\to\infty$. This may be written as
$\rho^w_p(f^*_n/\omega,1/w)\to 0$ as $n\to\infty$. 
Since $\log^+|f^*_n|\in L^p$, it follows that $f^*_n/\omega\in L^{\omega}_p$
for all $n$. Hence $f^*_n/\omega\to 1/w\notin L^{\omega}_p$ in the space
$(L^w_p,\rho^w_p)$, and therefore $L^{\omega}_p$ is not closed subspace
of the space $L^w_p$. This contradicts the assumption
that $(L^{\omega}_p,\rho^w_p)$ is complete. Hence must be 
$\log \left(w/{\omega}\right)\in L^p$, and the second assertion of
(iv) follows immediately from (i) and (ii).

\begin{cor}
If $w$ is a weight such that  $\log^+w\in L^p$, then  $L_p\subset L^w_p$, 
and $L_p= L^w_p$ if and only if $\log w\in L^p$. Furthermore,
 the topology  defined on $L_p$ by the metric $\rho_p$
 is stronger than that induced on $L_p$ by the metric $\rho^w_p$,
and if $\log w\in L^p$, then these two topologies coincide.
\end{cor}
{\it Proof.}
The proof follows directly from the assertions (i) and (ii) of Theorem 6.3,
by putting $\omega(t)\equiv 1$ on $[0,2\pi)$. 

\begin{thm} If $\log w\in L^p$, then  $(L_p,\rho^w_p)$ is an 
$F$-algebra, and the metric topologies  $\rho^w_p$ and $\rho_p$ are the same.
\end{thm}
{\it Proof.}
Since by Corollary 2.4, $(L_p,\rho_p)$ is an $F$-algebra, to show that 
$(L_p,\rho^w_p)$ is an $F$-algebra, it suffices by Corollary 6.4 to prove that 
$(L_p,\rho^w_p)$ is complete.
Let $(f_n)$ be a Cauchy sequence in $(L_p,\rho^w_p)$. This means
that $(wf_n)$ is a Cauchy sequence in $(L_p,\rho_p)$, and by Corollary 2.4,
there is a $g\in L_p$ such that $\rho_p(wf_n,g)\to 0$ as $n\to\infty$.
Then $\rho^w_p(f_n,g/w)=\rho_p(wf_n,g)\to 0$, and since $\rho^w_p(g/w,0)=
\rho_p(g,0)<\infty$, we see that $f=g/w\in L^w_p=L_p$. Hence, the space
$(L_p,\rho^w_p)$ is complete, and the theorem is proved.\\

\noindent{\it Remark.} {\it It is easy to see that the vector space $L^w_p$ is 
not necessarily algebra with respect to pointwise addition and multiplication.}

\begin{cor} The initial metric topology $\rho_p$
is a unique $\rho^{\omega}_p$ metric topology with $\log^+w\in L^p$
that makes $(L_p,\rho^w_p)$ into an $F$-space.
\end{cor}
{\it Proof.}
Suppose that $\log^+w\in L^p$, and that $(L_p,\rho^w_p)$ is an $F$-space.
 Consider the identity map
\[
j:(L_p,\rho_p)\to (L_p,\rho^w_p).
\] 
  Then by Corollary 6.4, $j$ is continuous. Since by Theorem 6.2
  $(L_p,\rho_p)$ is an $F$-space, by the open mapping theorem, we conclude that
 the inverse $j^{-1}$ of $j$ is also continuous. This shows that the metric 
topologies $\rho_p$ and $\rho^w_p$ must be the same on $L_p$.

\begin{thm} Let  $p\ge 1$ and let $w$ be a weight such that $\log^+w\in L^p$. 
Then the following statements about $w$ are equivalent.

{\rm (i)} $\log w\in L^p$;

{\rm (ii)} $L^w_p=L_p$;

{\rm (iii)} The metrics $\rho^w_p$ and $\rho_p$ define the same topology on 
$L_p$;

{\rm (iv)} $(L_p,\rho^w_p)$ is a complete metric space.

{\rm (v)} $(L_p,\rho^w_p)$ is an $F$-algebra;
\end{thm}
{\it Proof.}

(i)$\Leftrightarrow$(ii). Follows from Corollary 6.4.

(iii)$\Rightarrow$(i). This is a consequence of the assertion (iii)
of Theorem 6.3, by setting $\omega(t)\equiv 1$ for $t\in [0,2\pi)$.   

(i)$\Rightarrow$(iii). This is immediate from Corollary 6.4.

(iv)$\Rightarrow$(i). This follows from the assertion (iv)
of Theorem 6.3, by setting $\omega(t)\equiv 1$ for $t\in [0,2\pi)$.

(i)$\Rightarrow$(v). This is immediate  from Corollary 6.5.

(v)$\Rightarrow$(iv). This is obvious.

\vspace{5mm}

${}^{*}$  Maritime Faculty, University of Montenegro,
85330 Kotor, Montenegro

\quad {\it E-mail address}: romeo@ac.me

\vspace{2mm}

${}^{**}$  Faculty of Science, 
University of Montenegro, 81000 Podgorica, Montenegro 

\quad {\it E-mail address}: zarkop@ac.me

\vspace{2mm}

${}^{***}$ Faculty of Science, 
University of Montenegro, 81000 Podgorica, Montenegro 


\vspace{5mm}

\begin{thebibliography}{99}
\bibitem[DS]{1} N.~Dunford and J.~T.~Schwartz, \emph{Linear operators} I,
Wiley--Interscience, New York, 1958.  
 \bibitem[D]{2} P.~L.~Duren, \emph{Theory of $H^p$ spaces}, Academic Press, New
York, 1970.
\bibitem[G]{3}
T.~W.~Gamelin, \emph{Uniform algebras}, Prentice--Hall, Englewood Cliffs,
New Jersey, 1969.
\bibitem[L]{4}
R.~Le\'{s}niewicz, \emph{On linear functionals in Hardy--Orlicz spaces},
I, Studia Math. \textbf{46} (1973), 53--77. 
\bibitem[MP]{5} R.~Me\v{s}trovi\'c and \v{Z}.~Pavi\'cevi\'c, 
\emph{The logarithmic analogue of Szeg\"{o}'s theorem}, Acta Sci. Math. 
(Szeged) \textbf{64} (1998), 97--102. 
\bibitem[Mo]{6}  N.~Mochizuki, \emph{Algebras of holomorphic  functions between 
$H^p$ and $N_*$}, Proc. Amer. Math. Soc. \textbf{105} (1989), 898--902.
\bibitem[Mu]{7}
J.~Musielak, \emph{Orlicz spaces and modular spaces}, Lecture Notes
in Math. \textbf{1034}, Springer--Verlag, 1983.
\bibitem[MO]{8}
J.~Musielak and W.~Orlicz,  \emph{Some remarks on modular spaces}, Bull. Acad.
 Polon. Sci., \textbf{7} (1959), 661--668.
\bibitem[S]{9} M.~Stoll, \emph{ Mean growth and Taylor coefficients of some 
topological algebras of analytic functions}, Ann. Polon. Math. \textbf{35} 
(1977), 139--158. 
\bibitem[P]{10}
I.~I.~Privalov, \emph{Boundary properties of analytic functions},
Moscow University Press, Moscow, 1941; 2nd ed., GITTL, Moscow, 1950. (Russian)
\bibitem[Y]{11} N.~Yanagihara, \emph{Multipliers and linear functionals for the 
class $N^+$}, Trans. Amer. Math. Soc. \textbf{180} (1973), 449--461. 
\end{thebibliography}
\end{document}